%% file: SIAM.tex
\documentclass[hidelinks,onefignum,onetabnum]{siamart220329}

\usepackage{lipsum}
\usepackage{amsfonts}
\usepackage{graphicx}
\usepackage{epstopdf}

\usepackage{algorithm}
\usepackage{algorithmicx}
\usepackage[noend]{algpseudocode}

\ifpdf
  \DeclareGraphicsExtensions{.eps,.pdf,.png,.jpg}
\else
  \DeclareGraphicsExtensions{.eps}
\fi

\newsiamremark{remark}{Remark}
\newsiamremark{hypothesis}{Hypothesis}
\crefname{hypothesis}{Hypothesis}{Hypotheses}
\newsiamthm{claim}{Claim}

\headers{Faster randomized partial trace estimation}{T. Chen, R. Chen, K. Li, S. Nzeuton, Y. Pan, and Y. Wang}

\title{Faster randomized partial trace estimation\thanks{Submitted to the editors DATE.
\funding{This work was supported in part through the NYU IT High Performance Computing resources, services, and staff expertise.}}}

\author{
Tyler Chen\thanks{New York University, New York , NY (\email{tyler.chen@nyu.edu}, \url{https://research.chen.pw})}
\and Robert Chen\footnotemark[2]
\and Kevin Li\footnotemark[2]
\and Skai Nzeuton\thanks{Cornell University}
\and Yilu Pan\thanks{New York University, Shanghai}
\and Yixin Wang\footnotemark[2]  
}

\usepackage{amsopn}

\newcommand{\tr}{\operatorname{tr}}
\usepackage{bm}

\usepackage{dsfont}

\newcommand{\ii}{\mathrm{i}}

\newcommand{\EE}{\mathbb{E}}

\newcommand{\VV}{\mathbb{V}}

\renewcommand{\vec}{\mathbf}
\newcommand{\cT}{\dagger}
\newcommand{\T}{\mathsf{T}}
\newcommand{\F}{\mathsf{F}}

\newcommand{\syss}{\mathrm{s}}
\newcommand{\sysb}{\mathrm{b}}
\newcommand{\syst}{\mathrm{t}}

\newcommand{\ds}{d_{\syss}}
\newcommand{\db}{d_{\sysb}}
\newcommand{\dt}{d_{\syst}}

\newcommand{\Hs}{\mathcal{H}_{\syss}}
\newcommand{\Hb}{\mathcal{H}_{\sysb}}
\newcommand{\Ht}{\mathcal{H}_{\syst}}

\newcommand{\trb}{\tr_{\sysb}}
\newcommand{\trbest}{\widehat{\mathsf{tr}}_{\sysb}}

\usepackage{xcolor}

\definecolor{c0}{HTML}{140e36}
\definecolor{c1}{HTML}{641a80}
\definecolor{c2}{HTML}{f7705c}
\definecolor{c3}{HTML}{febd82}

\usepackage{tikz}

\usepackage{hyperref}
\hypersetup{
    colorlinks,
    linkcolor=c1,
    citecolor=c1,
    urlcolor=c1
}

\ifpdf
\hypersetup{
  pdftitle={Faster randomized partial trace estimation},
  pdfauthor={T. Chen, R. Chen, K. Li, S. Nzeuton, Y. Pan, and Y. Wang}
}
\fi

\begin{document}

\maketitle


\begin{abstract}
    We develop randomized matrix-free algorithms for estimating partial traces, a generalization of the trace arising in quantum physics and chemistry.
    Our algorithm improves on the typicality-based approach used in [T. Chen and Y-C. Cheng, \emph{Numerical computation of the equilibrium-reduced density matrix for strongly coupled open quantum systems},  J. Chem. Phys. 157, 064106 (2022)] by deflating important subspaces (e.g. corresponding to the low-energy eigenstates) explicitly.
    This results in a significant variance reduction, leading to several order-of-magnitude speedups over the previous state of the art.
    We then apply our algorithm to study the thermodynamics of several Heisenberg spin systems, particularly the entanglement spectrum and ergotropy.
\end{abstract}

\begin{keywords}
partial trace, deflation, stochastric trace, Krylov subspace method
\end{keywords}

\begin{MSCcodes}
68Q25, 68R10, 68U05
\end{MSCcodes}

\section{Introduction}


The state of a quantum system is described by a \emph{density matrix} of dimension exponential in the system size.
Often we are interested in the the state of a subsystem of the total system.
This can be obtained by taking the \emph{partial trace} of the total system density matrix, which yields a density matrix of dimension depending only on the size of the subsystem of interest (called the reduced system density matrix).
This matrix can then be used to understand important properties of the subsystem, for instance its entanglement with the rest of the total system \cite{nielsen_chuang_10}.
Recently, a number of numerical methods have been developed to estimate partial traces \cite{wang_ashhab_nori_11,chiu_strathearn_keeling_22,chen_cheng_22}.

If the total system density matrix is known explicitly, computing the partial trace is trivial.
Unfortunately, owing to the exponential dependence of the total system density matrix on the system size, it is typically prohibitively expensive to obtain (or even store) the total system density matrix.
Hope is not lost; in many situations, the total system density matrix has an implicit representation in terms of a (typically sparse) Hamiltonian $\vec{H}$ describing the configuration of the system of interest. 
For instance, the total system density matrix may be proportional to $\exp(-\beta \vec{H})$ for some parameter $\beta>0$, or might be obtained from $\vec{H}$ by solving the Schr\"odinger equation with a given initial condition.

Mathematically, this means that the task of computing reduced density matrices can often be viewed as computing the partial trace of a \emph{matrix function} of $\vec{H}$.
In many situations, the Hamiltonian is sparse and admits implicit matrix-vector products (i.e we can compute the map $\vec{x} \mapsto \vec{H}\vec{x}$).
Thus, for moderately sized systems for which it is possible to store dense vectors, Krylov subspace methods offer an attractive potential approach.
Indeed, such methods are widely used for the closely related task of trace estimation \cite{skilling_89,jaklic_prelovsek_94,schnalle_schnack_10,weisse_wellein_alvermann_fehske_06,schnack_richter_steinigeweg_20,schulter_gayk_schmidt_honecker_schnack_21,jin_willsch_willsch_lagemann_michielsen_deraedt_21,baer_neuhauser_rabani_22}.

\section{Background}
We begin by providing some background on a natural setting in which partial traces arise, as well as on existing methods for implicit partial trace approximation.

\subsection{Equilibrium reduced density matrices}
Throughout, the total system is defined on a finite dimensional Hilbert space $\Ht = \Hs\otimes \Hb$, where $\Hs$ and $\Hb$ are the Hilbert spaces for subsystem (s) and subsystem (b) respectively.
We assume the total system is governed by a Hamiltonian
\begin{equation}
\vec{H} = \bar{\vec{H}}_{\syss} + \bar{\vec{H}}_{\sysb} + \vec{H}_{\syss\sysb},
\end{equation}
where $\bar{\vec{H}}_{\syss} = \vec{H}_{\syss} \otimes  \vec{I}_{\sysb}$ corresponds to the Hamiltonian of subsystem (s), $\bar{\vec{H}}_{\sysb} = \vec{I}_{\syss} \otimes  \vec{H}_{\sysb}$ corresponds to the Hamiltonian of subsystem (b), and $\vec{H}_{\syss\sysb}$ accounts for \emph{non-negligible} interactions between the two subsystems.

When the total system is in thermal equilibrium at inverse temperature $\beta$ (due to weak coupling with a ``superbath''), the state of the system is described by a density matrix 
\begin{equation}
\label{eqn:rho_total}
    \bm{\rho}_{\syst} = \bm{\rho}_{\syst}(\beta) = \frac{\exp(-\beta \vec{H})}{Z_{\syst}}
    ,\qquad
    Z_{\syst} = Z_{\syst}(\beta) = \tr(\exp(-\beta \vec{H}));
\end{equation}
see for instance \cite{gemmer_michel_mahler_09,vinjanampathy_anders_16,alicki_kosloff_18}.
The quantity $Z_{\syst}(\beta)$ is called the partition function and provides insight into a number of thermodynamic properties of the system.

Often, we are interested in the state of subsystem (s) rather than the total system.
If subsystem (s) did not interact with subsystem (b) (i.e. if $\vec{H}_{\syss\sysb} = \vec{0}$), then the density matrix for subsystem (s) would simply be proportional to $\exp(-\beta \vec{H}_{\syss})$.
However, when the interactions between subsystems (s) and (b) are non-negligible, the density matrix $\bm{\rho}^*$ for subsystem (s) is instead obtained by ``tracing out''\footnote{This is analogous to ``integrating out'' variables from a joint probability distribution to obtain the marginal distribution for a variable of interest.} the effects of subsystem (b); i.e.
\begin{equation}
\label{eqn:rho_HMF}
\bm{\rho}^* 
= \bm{\rho}^* (\beta)
= \trb(\bm{\rho}_{\syst})
= \frac{\trb(\exp(-\beta \vec{H}))}{\tr(\exp(-\beta \vec{H}))},
\end{equation}
where $\tr_b(\:\cdot\:)$ is the \emph{partial trace} over subsystem (b) \cite{campisi_zueco_talkner_10,ingold_hanggi_talkner_09,talkner_hanggi_20}.

\subsection{Partial traces}
Let $\ds$ and $\db$ be the dimension of $\Hs$ and $\Hb$ respectively, so that $\dt = \ds\db$ is the dimension of $\Ht = \Hs\otimes\Hb$.
A general matrix $\vec{A} : \Ht \to \Ht$ can be partitioned as
\begin{align}
    \vec{A}
    = 
    \begin{bmatrix}
    \vec{A}_{1,1} & \vec{A}_{1,2} & \cdots & \vec{A}_{1,\ds} \\
    \vec{A}_{2,1} & \vec{A}_{2,2} & \cdots & \vec{A}_{2,\ds} \\
    \vdots & \vdots & \ddots & \vdots \\
    \vec{A}_{\ds,1} & \vec{A}_{\ds,2} & \cdots & \vec{A}_{\ds,\ds}
    \end{bmatrix},
\end{align}
where $\vec{A}_{i,j} : \Hb \to \Hb$ for each $i,j$.
The partial trace of $\vec{A}$ over $\Hb$ is defined as
\begin{equation}
    \label{eqn:partial_trace_def}
    \trb(\vec{A}) :=
    \begin{bmatrix}
    \tr(\vec{A}_{1,1}) & \tr(\vec{A}_{1,2}) & \cdots & \tr(\vec{A}_{1,\ds}) \\
    \tr(\vec{A}_{2,1}) & \tr(\vec{A}_{2,2}) & \cdots & \tr(\vec{A}_{2,\ds}) \\
    \vdots & \vdots & \ddots & \vdots \\
    \tr(\vec{A}_{\ds,1}) & \tr(\vec{A}_{\ds,2}) & \cdots & \tr(\vec{A}_{\ds,\ds})
    \end{bmatrix}.
\end{equation}
From \cref{eqn:partial_trace_def} it is clear that the partial trace is easy to compute if we have an explicit representation of $\vec{A}$.
However, we are most interested in the case $\vec{A} = \exp(-\beta \vec{H})$, where $\vec{H}$ is so large that storing and/or computing an explicit representation of $\vec{A}$ is intractable.
As such, we will consider only methods which access $\vec{A}$ through matrix-vector products which can then be approximated using Krylov subspace methods.

\subsection{Stochastic trace estimation}
Consider a real matrix $\vec{M} : \Hb\to\Hb$.
If $\vec{v} \in \Hb$ is a random vector whose entries are independent and identically distributed (iid) standard real Gaussian random variables, it is straightforward to show \cite{martinsson_tropp_20} that 
\begin{align}
\EE\bigl[\vec{v}^\T \vec{M} \vec{v}\bigr] = \tr(\vec{M})
,\qquad
\VV\bigl[\vec{v}^\T \vec{M} \vec{v} \bigr] 
= \frac{1}{2}\| \vec{M} + \vec{M}^\T \|_\F^2
\leq 2 \| \vec{M}\|_\F^2.
\end{align}
The use of estimators of this form\footnote{In quantum physics, such estimators are closely related to the idea of quantum typicality \cite{schrodinger_27,vonneumann_29}, which refers to the idea that, in many cases, a random state is representative of the overall state of a system; see \cite{goldstein_lebowitz_mastrodonato_tumulka_zanghi_10} for a review.} (although possibly with a different distribution) for approximating the trace of implicit matrices has been used since the late 1980s \cite{girard_87,skilling_89,hutchinson_89}.
Theoretical tail bounds appear in the physics \cite[etc.]{popescu_short_winter_06,reimann_07,gogolin_10} and numerical analysis \cite[etc.]{avron_toledo_11,roostakhorasani_ascher_14,meyer_musco_musco_woodruff_21,cortinovis_kressner_21} literature. 
These bounds control the probability that the estimator $\vec{v}^\T \vec{M}\vec{v}$ is far from the trace in terms of properties of $\vec{M}$ such as $\|\vec{M}\|_\F$ and $\|\vec{M}\|_2$.

\begin{figure}[t]
    \centering
    \includegraphics[scale=.55]{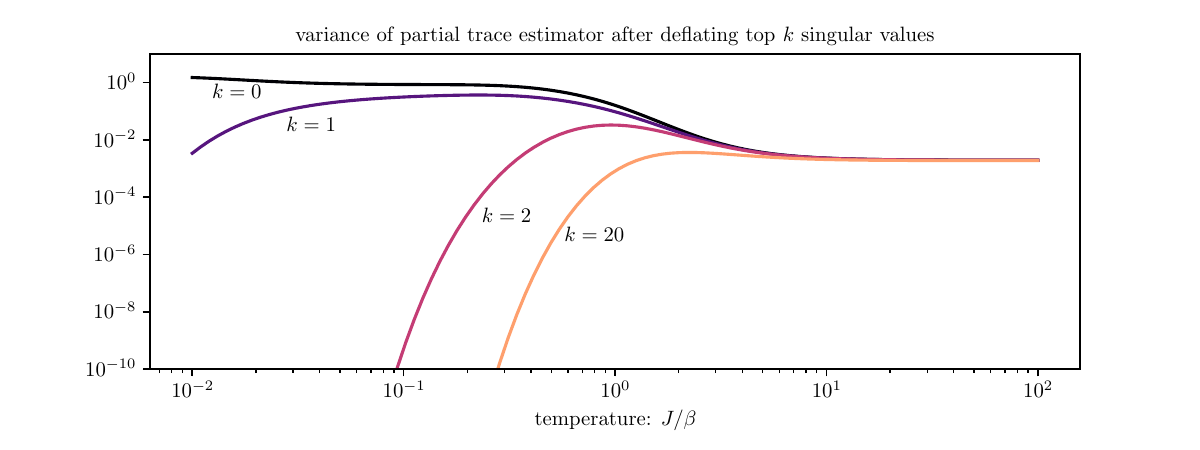}
    \caption{
    Let $\sigma_1 \geq \sigma_2 \geq \cdots \geq \sigma_{\dt}$ be the singular values of $\vec{A}$. 
    Then the variance of the partial trace estimator \cref{eqn:quadratic_partial_trace_estimator} is bounded above by
    $2\|\vec{A}\|_\F^2 = 2(\sigma_1^2 + \cdots + \sigma_{\dt}^2)$. 
    By deflating the top $k$ singular values, we can reduce the variance to at most $2(\sigma_{k+1}^2 + \cdots + \sigma_{\dt}^2)$ (see \cref{sec:variance_reduce}).
    Here we take $\vec{A} = \bm{\rho}_{\syst} = \exp(-\beta \vec{H}) / \tr(\exp(-\beta \vec{H}))$, where $\vec{H}$ corresponds to the solvable spin chain with $N=10$ and $h=0.3$ (described in \cref{sec:solvable_chain}).
    We then plot the variance bound $2(\sigma_{k+1}^2 + \cdots + \sigma_{\dt}^2)$ for several values of $k$ and a range of $\beta$.
    \emph{Takeaway}: The $k=0$ curve exhibits high variance at low-temperature. Through the use of deflation, the variance can be reduced significantly.
    }
    \label{fig:variance_example}
\end{figure}

Standard trace estimators can be extended to partial traces \cite{chen_cheng_22}. 
In particular, 
\begin{equation}
\label{eqn:quadratic_partial_trace_estimator}
  (\vec{I}_{\ds} \otimes \vec{v})^\T
    \vec{A}
    (\vec{I}_{\ds} \otimes \vec{v})
    =     
    \begin{bmatrix}
    \vec{v}^\T \vec{A}_{1,1} \vec{v} & \vec{v}^\T \vec{A}_{1,2} \vec{v} & \cdots & \vec{v}^\T \vec{A}_{1,\ds} \vec{v} \\
    \vec{v}^\T \vec{A}_{2,1} \vec{v} & \vec{v}^\T \vec{A}_{2,2} \vec{v} & \cdots & \vec{v}^\T \vec{A}_{2,\ds} \vec{v} \\
    \vdots & \vdots & \ddots & \vdots \\
    \vec{v}^\T \vec{A}_{\ds,1} \vec{v} & \vec{v}^\T \vec{A}_{\ds,2} \vec{v} & \cdots & \vec{v}^\T \vec{A}_{\ds,\ds} \vec{v}
    \end{bmatrix}
\end{equation}
provides an unbiased estimator for $\trb(\vec{A})$ when $\vec{v}$ is sampled as described above.
Given iid copies $\vec{v}_1, \ldots, \vec{v}_m$ of $\vec{v}$, we arrive at an estimator
\begin{equation}
\label{eqn:trbest}
    \trbest^{m}(\vec{A})
    := \frac{1}{m} 
    \sum_{i=1}^{m}
    (\vec{I}_{\ds} \otimes \vec{v}_i)^\T \vec{A} (\vec{I}_{\ds} \otimes \vec{v}_i).
\end{equation}
This estimator was studied in \cite{chen_cheng_22} for approximating reduced density matrices and will serve as the backbone of the algorithms developed in this paper.

The variance of a random matrix $\vec{X}$ can be defined as 
\begin{align}    
    \VV[\vec{X}] 
    := 
    \EE \Big[ \big\| \vec{X} - \EE[\vec{X}] \big\|_\F^2 \Big].
\end{align}
This is equivalent to the sum of the variances of the entries of $\vec{X}$, so assuming $\vec{v}$ has iid standard normal entries, we find that
\begin{equation}
\label{eqn:simple_est_var}
    \VV\Bigl[ \trbest^{m}(\vec{A}) \Bigr]
    = \frac{1}{m} \VV\Bigl[(\vec{I}_{\ds} \otimes \vec{v})^\T
    \vec{A}
    (\vec{I}_{\ds} \otimes \vec{v})\Bigr]
    \leq \frac{1}{m} \sum_{i=1}^{\ds}\sum_{j=1}^{\ds} 2\|\vec{A}_{i,j}\|_\F^2
    = \frac{2}{m} \| \vec{A} \|_\F^2.
\end{equation}

As with all Monte Carlo estimators which output a sample average, the estimator \cref{eqn:trbest} often suffers from large variance with fluctuations about the mean on the order of $\|\vec{A}\|_\F/\sqrt{m}$.
When computing quantities of the form $\trb(\vec{A})/\tr(\vec{A})$ these effects become particularly pronounced when $\vec{A}$ has a very quickly decaying spectrum.
Indeed, if $\vec{A}$ has one dominating (positive) eigenvalue, $\|\vec{A}\|_\F / \tr(\vec{A}) \approx 1$, whereas if $\vec{A}$ were close to the identity, then $\|\vec{A}\|_\F / \tr(\vec{A}) \approx 1/\sqrt{\dt}$.
In the case $\vec{A} = \bm{\rho}_{\syst} =  \exp(-\beta \vec{H}) / \tr(\exp(-\beta\vec{H}))$, this corresponds to difficulties when $\beta$ is very large (low temperature). 
We visualize how the variance depends on the temperature for this matrix in \cref{fig:variance_example}.

In the zero-temperature limit $\beta\to\infty$, the trace $\tr(\exp(-\beta \vec{H}))$ is determined entirely by the smallest (most negative) eigenvalue of $\vec{H}$, and the partial trace $\trb(\exp(-\beta \vec{H}))$ by the corresponding eigenvector (assuming a one dimensional eigenspace). 
This means a randomized estimator such as \cref{eqn:trbest} is not needed! 
Rather, one can simply apply classical techniques for obtaining extremal eigenvectors.
This paper is built on the fact that at low (but non-zero) temperatures, knowledge of the eigenvectors corresponding to small eigenvalues is still very useful.
In particular, we provide a deflation-based technique which can significantly reduce the variance of \cref{eqn:trbest} at low temperatures.
Our approach is closely related to \cite{morita_tohyama_20} and other deflation-based approaches for regular trace estimation \cite{weisse_wellein_alvermann_fehske_06,saibaba_alexanderian_ipsen_17,gambhir_stathopoulos_orginos_17}.
The potential for variance reduction is also visualized in \cref{fig:variance_example}.

\subsection{Contributions}
The primary contribution of this paper is a variance reduction technique for \cref{eqn:quadratic_partial_trace_estimator}.
This results in several order-of-magnitude speedups over the current state of the art algorithm for approximating partial traces of matrix functions \cite{chen_cheng_22}.
As such, we are able to study properties of quantum systems too large for existing methods.
While similar variance reduction approaches have been used for regular trace estimation \cite{girard_87,weisse_wellein_alvermann_fehske_06,lin_16,wu_laeuchli_kalantzis_stathopoulos_gallopoulos_16,gambhir_stathopoulos_orginos_17,morita_tohyama_20,meyer_musco_musco_woodruff_21}, partial traces do not satisfy a cyclic property which makes generalizing past work (numerically) difficult.
We propose some solutions to this difficulty.

Another contribution of this paper is to highlight some important problems from quantum physics which have been under-explored by the scientific computing community.
We feel that there is significant potential for increased collaboration between these communities which is limited by a lack cross-disciplinary knowledge transfer.

\section{A variance reduced algorithm}
\label{sec:variance_reduce}

We now describe a general technique for reducing the variance of the estimator \cref{eqn:trbest} for an arbitrary matrix $\vec{A}$.
By the linearity of the partial trace, for any matrix $\widetilde{\vec{A}}$,
\begin{equation}
\label{eqn:split}
    \trb(\vec{A})
    = \trb(\widetilde{\vec{A}}) + \trb(\vec{A} - \widetilde{\vec{A}}).
\end{equation}
If we are able to compute the partial trace of the first term exactly, we can estimate the partial trace of $\trb(\vec{A})$ by applying the randomized estimator \cref{eqn:quadratic_partial_trace_estimator} to the residual term. 
That is, by
\begin{equation}
\label{eqn:residual}
    \trb(\vec{A}) 
    \approx
    \trb(\widetilde{\vec{A}}) + \trbest^m(\vec{A} - \widetilde{\vec{A}}).
\end{equation}
The variance of such an estimate is entirely due to the variance of $\trbest^m(\vec{A} - \widetilde{\vec{A}})$.
Thus, if $\|\vec{A} - \widetilde{\vec{A}}\|_\F^2 < \| \vec{A} \|_\F^2$, then the variance of the estimator on the right of \cref{eqn:residual} is reduced over that of $\trbest^m(\vec{A})$.

Splittings similar to \cref{eqn:residual} have previously been used as a variance reduction technique for regular trace estimation
\cite{girard_87,weisse_wellein_alvermann_fehske_06,lin_16,wu_laeuchli_kalantzis_stathopoulos_gallopoulos_16,gambhir_stathopoulos_orginos_17,morita_tohyama_20}.
Perhaps the most widely known approach in the numerical analysis and theoretical computer science communities is the Hutch++ algorithm \cite{meyer_musco_musco_woodruff_21} which produces a $1 \pm \epsilon$ relative approximation to the trace of a positive semi-definite matrix using just $O(\epsilon^{-1})$ matrix-vector products with $\vec{A}$.
Several improvements to this algorithm have been proposed \cite{persson_cortinovis_kressner_22,epperly_tropp_webber_23}, including for the case $\vec{A} = f(\vec{H})$ \cite{persson_kressner_23,chen_hallman_23}.

\subsection{Partial trace of low-rank matrices}

The partial trace of rank-1 matrix can be computed efficiently given a factorization as an outer product.
In particular, for any $\vec{x} \in \Ht$ we can write the outer product as
\begin{equation}
        \vec{x}\vec{x}^\T = \begin{bmatrix}
        \vec{x}_{(1)}\vec{x}_{(1)}^\T & \vec{x}_{(1)}\vec{x}_{(2)}^\T & \cdots & \vec{x}_{(1)}\vec{x}_{(\ds)}^\T \\
        \vec{x}_{(2)}\vec{x}_{(1)}^\T & \vec{x}_{(2)}\vec{x}_{(2)}^\T & \cdots & \vec{x}_{(2)}\vec{x}_{(\ds)}^\T \\
        \vdots & \vdots & \ddots & \vdots \\
        \vec{x}_{(n_v)}\vec{x}_{(1)}^\T & \vec{x}_{(\ds)}\vec{x}_{(2)}^\T & \cdots & \vec{x}_{(\ds)}\vec{x}_{(\ds)}^\T 
    \end{bmatrix},
    \qquad
    \vec{x}
    = \begin{bmatrix}
        \vec{x}_{(1)} \\ \vec{x}_{(2)} \\ \vdots \\ \vec{x}_{(\ds)}
    \end{bmatrix}.
\end{equation}
Using that $\tr(\vec{x}_{(i)}\vec{x}_{(j)}^\T) = \vec{x}_{(i)}^\T \vec{x}_{(j)}$, we find that
\begin{equation}
\label{eqn:rank1_partial_trace}
    \trb(\vec{x}\vec{x}^\T) 
    = \begin{bmatrix}
        \vec{x}_{(1)}^\T \vec{x}_{(1)} & \vec{x}_{(2)}^\T \vec{x}_{(1)} & \cdots & \vec{x}_{(\ds)}^\T \vec{x}_{(1)} \\
        \vec{x}_{(1)}^\T \vec{x}_{(2)} & \vec{x}_{(2)}^\T \vec{x}_{(2)} & \cdots & \vec{x}_{(\ds)}^\T \vec{x}_{(2)} \\
        \vdots & \vdots & \ddots & \vdots \\
        \vec{x}_{(1)}^\T \vec{x}_{(\ds)} & \vec{x}_{(2)}^\T \vec{x}_{(\ds)} & \cdots & \vec{x}_{(\ds)}^\T \vec{x}_{(\ds)} 
    \end{bmatrix}.
\end{equation}
This observation and the linearity of the partial trace allow us to efficiently compute the partial trace of a generic symmetric rank-$k$ matrix
\begin{equation}
\label{eqn:low_rank_form}
    \widetilde{\vec{A}}
    = \sum_{i=1}^{k} \theta_i \vec{x}_i\vec{x}_i^\T
\end{equation}
given access to the factors $(\theta_i,\vec{x}_i)$, $i=1,2,\ldots, k$.

\subsection{Implicit partial trace estimation}
It is clear that choosing $\widetilde{\vec{A}}$ as a rank-$k$ approximation to $\vec{A}$ will suit our needs.
In particular, let $\vec{Q} \in \mathbb{R}^{d\times k}$ be a matrix with orthonormal columns ($\vec{Q}^\T \vec{Q} = \vec{I}_k)$ and define
\begin{equation}
\label{eqn:QQAQQ}
    \widetilde{\vec{A}} := \vec{Q} \vec{Q}^\T \vec{A} \vec{Q} \vec{Q}^\T.
\end{equation}
We can efficiently obtain a factorization of the form \cref{eqn:low_rank_form} using just $k$ matrix-vector products with $\vec{A}$. 
Indeed, form $\vec{G} := \vec{Q}^\T \vec{A} \vec{Q}$, compute an eigendecomposition $\vec{G} =: \sum_{i=1}^{k} \theta_i \vec{s}_i \vec{s}_i^\T$, and set $\vec{x}_i := \vec{Q} \vec{s}_i$ to obtain a rank-$k$ approximation to \cref{eqn:QQAQQ} of the form \cref{eqn:low_rank_form}.

Then, rewriting \cref{eqn:residual}, we arrive at a computationally feasible estimator,
\begin{equation}
\label{eqn:residual_low-rank}
    \trbest^{m}(\vec{A};\vec{Q}) :=
    \sum_{i=1}^{k} \theta_i \trb(\vec{x}_i\vec{x}_i^\T) + \trbest^m(\vec{A} - \vec{Q} \vec{Q}^\T \vec{A} \vec{Q} \vec{Q}^\T).
\end{equation}

We provide pseudocode for computing \cref{eqn:residual_low-rank} and a corresponding error estimate in \cref{alg:main}.
The total number of matrix-vector products with $\vec{A}$ is $k+m\ds$; $k$ products are used to compute $\vec{A}\vec{Q}$ in \cref{alg:main:QAQ}, and in each of the $m$ loops, $\ds$ products are used to compute $\vec{A}\vec{Y}$ in \cref{alg:main:Brem}.

\begin{algorithm}[ht]
\caption{Variance reduced partial trace estimation}\label{alg:main}
\fontsize{10}{14}\selectfont
\begin{algorithmic}[1]
\Procedure{partial-trace}{$\vec{A}, \vec{Q}, m$}
\State $\vec{\Theta},\vec{S}  = \textsc{eig}(\vec{Q}^\T \vec{A} \vec{Q})$ \label{alg:main:eig}
\label{alg:main:QAQ} \Comment{$k\times k$ matrix}
\State $\vec{X} = \vec{Q} \vec{S}$ \Comment{$\vec{X} = [\vec{x}_1, \ldots, \vec{x}_k]$}
\For{$i=1,2,\ldots, k$}
\State $\vec{B}_{\textup{defl}}^{(i)} = \theta_i \trb(\vec{x}_i\vec{x}_i^\T)$ \Comment{using \cref{eqn:rank1_partial_trace}}
\EndFor
\For{$i=1,2,\ldots, m$}
\State $\vec{Y} = (\vec{I}_{\ds} \otimes \vec{v})$, $\vec{v}$ is length $\db$ iid Gaussian vector
\State $\vec{B}_{\textup{rem}}^{(i)} = \vec{Y}^\T \vec{A} \vec{Y} - \vec{Y}^\T \vec{X} \vec{\Theta} \vec{X}^\T \vec{Y}$
\label{alg:main:Brem}
\EndFor
\State \Return $\trbest^{m}(\vec{A};\vec{Q}) = \sum_{i=1}^{r} \vec{B}_{\textup{defl}}^{(i)} + \frac{1}{m} \sum_{i=1}^{m}\vec{B}_{\textup{rem}}^{(i)}$
\EndProcedure
\end{algorithmic}
\end{algorithm}

\subsection{Choosing the projection space}

\cref{alg:main} takes as input the matrix $\vec{Q}$ which determines the projection space used for deflation, and the quality of the output depends strongly on $\vec{Q}$.
For regular trace estimation, a natural approach is to obtain $\vec{Q}$ by sketching $\vec{A}$ \cite{meyer_musco_musco_woodruff_21}.
While sketching can be used to generate $\vec{Q}$ for use in \cref{alg:main}, there are some practical difficulties for the case $\vec{A} =f(\vec{H})$, which is the focus of this paper. 
We discuss these difficulties in \cref{sec:partial_trace_function}.

Another reasonable choice is to take $\vec{Q}$ aligned with eigenvectors of $\vec{A}$.
In fact, the choice of $\vec{Q}$ with $k$ columns which minimizes the variance of \cref{alg:main} is to take $\vec{Q}$ as the eigenvectors corresponding to the $k$ eigenvalues of $\vec{A}$ with largest magnitude (i.e. corresponding to the top $k$ singular values of $\vec{A}$).
In this case
\begin{align}
    \| \vec{A} - \vec{Q}\vec{Q}^\T \vec{A} \vec{Q} \vec{Q}^\T \|_\F^2
    = \min_{\operatorname{rank}(\widetilde{\vec{A}}) = k}\| \vec{A} - \widetilde{\vec{A}} \|_\F^2
    = 
    \sum_{i=k+1}^{\dt} \sigma_i^2,
\end{align}
where $\{\sigma_i\}$ are the singular values of $\vec{A}$ arranged in non-increasing order.
In the case $\vec{A} = \exp(-\beta \vec{H})$, the singular values of $\vec{A}$ are $\exp(-\beta \lambda_{i})$, where $\lambda_i$ are the eigenvalues of $\vec{H}$. 
When $\beta$ is large (low-temperature), several of these singular values are significantly larger than the others, and deflation is effective at decreasing the norm.  

To illustrate this idea quantitatively, suppose that, for some fixed constants $c,\alpha\in(0,1)$ and $c',c''>0$,
\begin{equation}
c' (c\alpha )^{i} <\sigma_i < c''\, \alpha^i,\qquad i=1,2,\ldots, \dt.
\end{equation}
Then, if $\vec{Q}$ contains the $k$ eigenvectors of $\vec{A}$ corresponding to the largest magnitude eigenvalues, the variance reduced estimator \cref{eqn:residual_low-rank} satisfies
\begin{equation}
\label{eqn:exponential_deflation}
\VV\Big[ \trbest^{m}(\vec{A};\vec{Q}) \Big]^{1/2}
< C \, \alpha^{k}\, \VV\Big[ \trbest^{m}(\vec{A}) \Big]^{1/2},
\end{equation}
for some $C = C(\alpha,c,c',c'')$ that does not depend on $k$ or the dimension $\dt$.\footnote{
Since $n\geq 1$ and $c,\alpha\in(0,1)$, $\alpha^{2n}>0$ and $(c\alpha)^{2n} < (c\alpha)^2 < \alpha^2$.
Therefore,
\begin{equation}
\frac{\sum_{i=k+1}^{\dt} \sigma_i^2}{\sum_{i=1}^{\dt}\sigma_i^2}
< \frac{c''}{c'}\frac{\sum_{i=k+1}^{n} \alpha^{2i}}{\sum_{i=1}^{\dt} (c \alpha)^{2i}}
= \frac{c''}{c'}\frac{(\alpha^{2k} - \alpha^{2\dt})(1-(c\alpha)^2)}{c^2(1-\alpha^2)(1-(c \alpha)^{2\dt})}
< \frac{c''}{c'}\frac{\alpha^{2k} }{c^2(1-\alpha^2)^2}.
\end{equation}
Rearranging and taking a square root gives \cref{eqn:exponential_deflation}.
}
In other words, if the singular values of $\vec{A}$ decay exponentially, deflating the top $k$ eigenvalues results in an exponential decrease in the magnitude of the fluctuations of the partial trace estimator \cref{eqn:residual_low-rank} over the basic estimator \cref{eqn:trbest} from \cite{chen_cheng_22}.

\section{Partial traces of matrix functions}\label{sec:partial_trace_function}

The primary focus of this paper is on estimating $\trb(\exp(-\beta\vec{H}))$; i.e. the partial trace of a matrix proportional to the density matrix $\bm{\rho}$ describing the state of the total system in thermal equilibrium at inverse temperature $\beta$.
In this section, we describe an implementation of \cref{alg:main} for general $f(\vec{H})$, with a particular focus on the case $f(x) = \exp(-\beta x)$. 
In the case of the standard trace $\tr(f(\vec{H}))$, this is commonly addressed using a combination of the typicality estimator $\vec{v}^\T f(\vec{H})\vec{v}$ and Krylov subspace methods \cite{skilling_89,jaklic_prelovsek_94,weisse_wellein_alvermann_fehske_06,ubaru_chen_saad_17,han_malioutov_avron_shin_17,schnack_richter_steinigeweg_20,schulter_gayk_schmidt_honecker_schnack_21,jin_willsch_willsch_lagemann_michielsen_deraedt_21,chen_trogdon_ubaru_22}.

\Cref{alg:main} requires computing quantities like
\begin{equation}
\label{eqn:diff_hard}
    \vec{Y}^\T \vec{A} \vec{Y} - \vec{Y}^\T (\vec{Q}\vec{Q}^\T \vec{A}\vec{Q}\vec{Q}^\T)\vec{Y}, \qquad \vec{Y} := (\vec{I}_{\ds}\otimes \vec{v}),
\end{equation}
which can be difficult to compute accurately due to the potential for cancellation errors.
In particular, each term of the difference may be much larger than the difference, so an accurate approximation to each term in a relative sense need not yield a good relative (or even additive) approximation to the difference.
For \cref{alg:main}, where products with $\vec{A}$ are assumed to be exact, this is not an issue. 
However, efficient methods for computing products with $\vec{A} = \exp(-\beta \vec{H})$, such as time-stepping and Krylov subspace methods \cite{druskin_knizhnerman_89,gallopoulos_saad_92,moler_vanloan_03}, result in some level of approximation error.

We will consider mainly the case that $\vec{Q}$ contains eigenvectors of $\vec{H}$ and hence of $\vec{A} = f(\vec{H})$.
In \cref{sec:general_Q} we discuss how one may be able to avoid the cancellation errors of \cref{eqn:diff_hard} for other $\vec{Q}$ with orthonormal columns.
We also discuss some pros-and-cons of various choices of $\vec{Q}$ and potential implementation difficulties, particularly in the context of $f(x) = \exp(-\beta x)$.

Write the eigendecomposition of $\vec{H}$ as
\begin{equation} \label{eqn:H_split}
    \vec{H} = 
    \begin{bmatrix}
        \vec{Q} & \widehat{\vec{Q}}
    \end{bmatrix}
    \begin{bmatrix}
        \vec{\Lambda} \\ & \widehat{\vec{\Lambda}}
    \end{bmatrix}
    \begin{bmatrix}
        ~\vec{Q}^\T~ \\  \widehat{\vec{Q}}^\T
    \end{bmatrix},
\end{equation}
where $\vec{Q}$ contains $r$ eigenvectors and $\vec{\Lambda}$ the corresponding $r$ eigenvalues.
We will assume that $\vec{Q}$ and $\vec{\Lambda}$ can be computed exactly.
Since there are many black-box and problem dependent techniques for this, we do not discuss particular methods for obtaining these quantities.

Using \cref{eqn:H_split}, we see $f(\vec{H})$ can be decomposed as
\begin{equation}\label{eqn:fH_split}
f(\vec{H})
    = \vec{Q} \vec{Q}^\T f(\vec{\Lambda}) \vec{Q}\vec{Q}^\T +  \widehat{\vec{Q}} \widehat{\vec{Q}}^\T f(\widehat{\vec{\Lambda}}) \widehat{\vec{Q}}  \widehat{\vec{Q}}^\T.
\end{equation}
This implies that
\begin{equation}\label{eqn:QQfHQQ}
\vec{Q}\vec{Q}^\T f(\vec{H}) \vec{Q} \vec{Q}^\T
= \vec{Q} f(\vec{\Lambda}) \vec{Q}^\T,
\end{equation}
which can be easily computed given $\vec{Q}$ and $\vec{\Lambda}$.
Thus, using \cref{eqn:fH_split} and the fact that the orthogonal projector $\widehat{\vec{Q}}  \widehat{\vec{Q}}^\T$ can equivalently be written $(\vec{I} - \vec{Q}\vec{Q}^\T)$, we obtain an alternate expression of the difference \cref{eqn:diff_hard}:
\begin{equation}\label{eqn:ZfHZ}
\vec{Y}^\T  f(\vec{H})  \vec{Y} - \vec{Y}^\T  \vec{Q} \vec{Q}^\T f(\vec{H}) \vec{Q}\vec{Q}^\T  \vec{Y}
= \vec{Z}^\T f(\vec{H}) \vec{Z}
,\qquad \vec{Z} := (\vec{I} - \vec{Q}\vec{Q}^\T) \vec{Y}.
\end{equation}
This allows us to avoid cancellation errors by ensuring that our approximation to $\vec{Z}^\T f(\vec{H}) \vec{Z}$ respects the fact that $\vec{Z}$ is orthogonal to $\vec{Q}$.

\subsection{The Lanczos algorithm with deflation}

In order to approximate the quantity 
\begin{equation}
\vec{Z}^\T  f(\vec{H})  \vec{Z},
\qquad \vec{Z} := (\vec{I} - \vec{Q}\vec{Q}^\T) \vec{Y}
\end{equation}
we make use of the block-Lanczos algorithm with explicit deflation.
The block-Lanczos method implicitly constructs an orthonormal basis $\vec{V} = [\vec{V}_0, \ldots, \vec{V}_{t-1}]$ for the block-Krylov subspace
\begin{equation}
\operatorname{span}\{ \vec{Z}, \vec{H} \vec{Z}, \ldots, \vec{H}^{t-1} \vec{Z} \}
\end{equation}
such that for all $j=0,1,\ldots, t-1$,
\begin{equation}
\operatorname{span}\{\vec{V}_0, \ldots, \vec{V}_j \} 
= \operatorname{span}\{ \vec{Z}, \vec{H} \vec{Z}, \ldots, \vec{H}^{j} \vec{Z} \}.
\end{equation}
In addition, the algorithm outputs a symmetric block-tridiagonal matrix 
\begin{equation}\label{eqn:Tdef}
    \vec{T} = \begin{bmatrix}
    \vec{M}_0 & \vec{R}_1^\T \\
    \vec{R}_1 & \ddots & \ddots \\
    & \ddots & \ddots & \vec{R}_{t-2}^\T \\
    &&\vec{R}_{t-2} & \vec{M}_{t-1}
    \end{bmatrix}
\end{equation}
satisfying $\vec{T} = \vec{V}^\T \vec{H}\vec{V}$.
Here we assume that the block-Krylov subspace does not become degenerate.

Since the input $\vec{Z}$ is orthogonal to the eigenvectors $\vec{Q}$, in exact arithmetic $\vec{V}$ will be entirely orthogonal to $\vec{Q}$  as well.
However, in finite precision arithmetic this cannot be guaranteed, and rounding errors might introduce small components in the directions of the $\vec{Q}$.
These errors can grow rapidly.
Thus, the block-Lanczos algorithm should be implemented to explicitly maintain orthogonality against $\vec{Q}$.
Such an implementation is given in \cref{alg:block_lanczos}.

\begin{algorithm}
\caption{Block-Lanczos algorithm with deflation}
\label{alg:block_lanczos}
\fontsize{10}{10}\selectfont
\begin{algorithmic}[1]
\Procedure{block-Lanczos-defl}{$\vec{H},\vec{Z},\vec{Q},t$}
\State \( \vec{V}_0,\vec{R}_0 = \textsc{qr}(\vec{Z}-\vec{Q}\vec{Q}^\T \vec{Z}) \),
\For {\( j=0,1,\ldots,t-1 \)}
    \State \( \vec{X} = \vec{H} \vec{V}_{j} - \vec{V}_{j-1} \vec{R}_{j-1}^\T \) \Comment{if $j=0$, $\vec{X} = \vec{H} \vec{V}_0$}
    \State \( \vec{M}_j = \vec{V}_{j}^\T \vec{X}  \)
    \State \( \vec{X} = \vec{X} - \vec{V}_{j} \vec{M}_j \)
    \State \( \vec{X} = \vec{X} - \vec{Q} \vec{Q}^\T \vec{X} \) \Comment{Explicit deflation}
    \State \label{alg:reorth}optionally, reorthogonalize $\vec{X}$ against $\vec{V}_0, \ldots, \vec{V}_{j-1}$
    \State \( \vec{V}_{j+1},\vec{R}_{j} = \textsc{qr}(\vec{X}) \)
\EndFor
\State \Return $\vec{T},\vec{R}_0$ \Comment{$\vec{T}$ as defined in \cref{eqn:Tdef}}
\EndProcedure
\end{algorithmic}
\end{algorithm}

In exact arithmetic, the Lanczos approximation to $\vec{Z}^\T  f(\vec{H})  \vec{Z}$ is given by
\begin{equation}\label{eqn:ZfHZ_lanczos}
\vec{Z}^\T  f(\vec{H})  \vec{Z} \approx \vec{R}_0^\T \vec{E}_1^\T f(\vec{T}) \vec{E}_1 \vec{R}_0,
\end{equation}
where $\vec{E}_1 = \vec{e}_1 \otimes \vec{I}$ and $\vec{R}_0$ is the R factor in the QR factorization of $\vec{Z}$.
This approximation is a block-Gauss quadrature approximation and is exact if $f$ is polynomial of degree at most $2t-1$ \cite[\S 6.6]{golub_meurant_09} \cite[Theorem 2.7]{frommer_lund_szyld_20}.
From this, we can obtain a simple bound for the convergence.

Let $\mathsf{Error} = \| \vec{Z}^\T  f(\vec{H})  \vec{Z}  - \vec{R}_0^\T \vec{E}_1^\T f(\vec{T}) \vec{E}_1 \vec{R}_0 \|$ and suppose $p$ is a polynomial of degree at most $2t-1$.
Then, with $\vec{P}_{\vec{Q}} = \vec{I} - \vec{Q}\vec{Q}^\T$,
\begin{align}
\mathsf{Error} & = \| \vec{Z}^\T  f(\vec{H})  \vec{Z} - \vec{Z}^\T  p(\vec{H})  \vec{Z} + \vec{R}_0^\T \vec{E}_1^\T p(\vec{T}) \vec{E}_1 \vec{R}_0  - \vec{R}_0^\T \vec{E}_1^\T f(\vec{T}) \vec{E}_1 \vec{R}_0 \|
\\& \leq \| \vec{Z}^\T  f(\vec{H})  \vec{Z} - \vec{Z}^\T  p(\vec{H})  \vec{Z} \| + \|  \vec{R}_0^\T \vec{E}_1^\T p(\vec{T}) \vec{E}_1 \vec{R}_0  - \vec{R}_0^\T \vec{E}_1^\T f(\vec{T}) \vec{E}_1 \vec{R}_0 \|
\\&\leq  \| \vec{Z} \|^2 \| \vec{P}_{\vec{Q}}( f(\vec{H}) - p(\vec{H}))\vec{P}_{\vec{Q}} \| + \| \vec{R}_0 \|^2 \| f(\vec{T}) - p(\vec{T}) \|.
\end{align}
Note that, with $\hat{\Lambda}$ denoting the diagonal entries of $\hat{\vec{\Lambda}}$ defined in \cref{eqn:H_split},
\begin{equation}
\| \vec{P}_{\vec{Q}}( f(\vec{H}) - p(\vec{H}))\vec{P}_{\vec{Q}} \| 
= \max_{x \in {\hat{\Lambda}}}| f(x) - p(x)|.
\end{equation}
Likewise, since $\vec{T} = \vec{V}^\T \vec{H} \vec{V} = \vec{V}^\T \vec{P}_{\vec{Q}}\vec{H}\vec{P}_{\vec{Q}} \vec{V} $, the eigenvalues of $\vec{T}$ are contained in the convex closure $\operatorname{conv}(\hat{\Lambda})$ of $\hat{\Lambda}$.
Thus,
\begin{equation}
\| f(\vec{T}) - p(\vec{T}) \|
= \max_{x \in \operatorname{spec}(\vec{T})}| f(x) - p(x)|
\leq \max_{x \in \operatorname{conv}(\hat{\Lambda})}| f(x) - p(x)|.
\end{equation}

Then, using that $\| \vec{Z} \| = \| \vec{R}_0 \|$ and that $p$ was arbitrary, we obtain the bound
\begin{equation}
\label{eqn:ZfHZ_bound}
\| \vec{Z}^\T  f(\vec{H})  \vec{Z}  - \vec{R}_0^\T \vec{E}_1^\T f(\vec{T}) \vec{E}_1 \vec{R}_0 \| 
\leq 2\| \vec{Z} \| \min_{\deg(p)<2t} \max_{x\in\operatorname{conv}(\hat{\Lambda})} | f(x) - p(x)|.
\end{equation}
Without deflation, note that $\hat{\Lambda} = \Lambda$, the set of eigenvalues of $\vec{H}$.
However, even when $r$ is small, $\operatorname{conv}(\hat{\Lambda})$ can be much smaller than $\operatorname{conv}(\Lambda)$. 
In such cases, deflation not only helps with variance reduction of the partial-trace estimator but also with the matrix function approximation.
This will also be an important consideration when we discuss the use of other projection spaces in \cref{sec:general_Q}.

\subsubsection{A note on finite precision arithmetic}
In exact arithmetic, the reorthgonalization step of \cref{alg:block_lanczos} is unnecessary as $\vec{X}$ is already orthogonal to $\vec{V}_0, \ldots, \vec{V}_{j-1}$. However, in finite precision arithmetic, failure to orthogonalize against these vectors at each iteration can lead to a drastic loss of orthogonality in $\vec{V}$.
In practice the formal expression \cref{eqn:ZfHZ} still converges in all examples we have observed.
This has been rigorously justified for the standard Lanczos method for approximating quadratic forms of matrix functions $(\ds=1)$ without deflation \cite{knizhnerman_96}.
The analysis in \cite{knizhnerman_96} is based on a careful analysis of the Lanczos algorithm in finite precision arithmetic \cite{paige_76,paige_80}. 
However, to the best of our knowledge, there is no similar analysis of the block Lanczos algorithm or the Lanczos algorithm with deflation.

\subsection{Algorithm}

We now have all the tools required to implement a version of \cref{alg:main} in the case $\vec{A} = f(\vec{H})$.
The resulting algorithm is summarized in \cref{alg:fH}.
In the context of regular trace estimation ($\ds=1$), similar approaches have been used successfully \cite{gambhir_stathopoulos_orginos_17,weisse_wellein_alvermann_fehske_06,morita_tohyama_20}. 
In particular, \cite{morita_tohyama_20} studies the task of of computing the partition function $Z(\beta)$ for a range of $\beta$, a task closely related to our main application of focus.
\begin{algorithm}[ht]
\caption{Variance reduced partial trace estimation for matrix functions}\label{alg:fH}
\fontsize{10}{14}\selectfont
\begin{algorithmic}[1]
\Procedure{partial-trace-func}{$\vec{H},f,k,m,t$}
\State Compute eigenvectors/values $\vec{Q},\vec{\Lambda}$ \label{alg:fH:eig}
\Comment{$\vec{Q} = [\vec{q}_1, \ldots, \vec{q}_k]$, $\vec{\Lambda} = \operatorname{diag}(\lambda_1, \ldots, \lambda_k)$}
\For{$i=1,2,\ldots, k$}
\State $\vec{B}_{\textup{defl}}^{(i)} = f(\lambda_i) \trb(\vec{q}_i\vec{q}_i^\T)$ \Comment{using \cref{eqn:rank1_partial_trace}}
\EndFor
\For{$i=1,2,\ldots, m$}\label{alg:fH:loop}
\State $\vec{Y} = (\vec{I}_{\ds} \otimes \vec{v})$, $\vec{v}$ is length $\db$ iid Gaussian vector
\State $\vec{Z} = (\vec{I} - \vec{Q}\vec{Q}^\T) \vec{Y}$ \Comment{Deflation}
\State $\vec{T},\vec{R}_0 = \textsc{block-Lanczos-defl}(\vec{H},\vec{Z},\vec{Q},t)$ \label{alg:fH:bld}
\State $\vec{B}_{\textup{rem}}^{(i)} = \vec{R}_0^\T \vec{E}_1^\T f(\vec{T}) \vec{E}_1 \vec{R}_0$
\Comment{$\vec{E}_1 = \vec{e}_1 \otimes \vec{I}$}
\EndFor
\State \Return $\trbest^{m}(\vec{A};\vec{Q}) = \sum_{i=1}^{k} \vec{B}_{\textup{defl}}^{(i)} + \frac{1}{m} \sum_{i=1}^{m}\vec{B}_{\textup{rem}}^{(i)}$
\EndProcedure
\end{algorithmic}
\end{algorithm}

\subsubsection{Computational costs}

\Cref{alg:fH} requires computing the $k$ eigenvalues/vectors of $\vec{H}$, the cost of which is context dependent.
The remaining number of matrix-vector products with $\vec{H}$ is $mtd_s$: in each of the $m$ loops, $d_s$ products are required to compute $\mathbf{H}\mathbf{V}_j$ in each of the $t$ iterations of the Block-Lanczos algorithm. 
The parameters $m$ and $t$ respectively control the the statistical variance of the partial trace estimator and the accuracy with which products matrix functions are computed. 
We note that the matrix-vector products for each of the $m$ samples can be done in parallel.

\subsubsection{Limitations and extensions}

The Lanczos-based method described in this section requires the storage of roughly $\ds$ dense vectors of length $\dt$, as well as repeated matrix-vector products with the total system Hamiltonian $\vec{H}_{\syst}$.
Since $\dt$ depends exponentially on the system size, this approach is only viable for moderately sized systems far from the thermodynamic limit.

The partial trace estimator \cref{eqn:quadratic_partial_trace_estimator} uses $\vec{v}\in\mathcal{H}_{\sysb}$ drawn from the uniform distribution on the hypersphere of radius $\sqrt{\db}$. 
However, the analogous estimator is still unbiased so long as  $\EE[\vec{v}\vec{v}^\T] = \vec{I}_{\sysb}$.
This opens the possibility of using an appropriate distribution on tensor network states \cite{orus_14,orus_19}.
While tensor network versions of the Lanczos algorithm have been studied \cite{dargel_wollert_honecker_mcculoch_schollwock_pruschke_12}, a more common imaginary time evolution approach \cite{orus_14,phien_mcculoch_vidal_15,czarnik_rams_dziarmaga_16,kshetrimayum_rizzi_eisert_orus_19} is likely suitable for approximating the action of $\exp(-\beta \vec{H}_{\syst})$ for sufficiently low-temperatures, at least as long as the total system has sufficiently local interactions.

\subsubsection{Using arbitrary projection matrices}
\label{sec:general_Q}

We may hope that we can obtain a matrix $\vec{Q}$ which reduces the variance of the partial trace estimator nearly as much as using the exact top eigenspace more efficiently than computing the top eigenspace exactly.
For any matrix $\vec{Q}$ with orthonormal columns it can be verified that
\begin{equation}
\vec{A} - \vec{Q}\vec{Q}^\T \vec{A}\vec{Q}\vec{Q}^\T
= \frac{1}{2} \Big[ (\vec{I} + \vec{Q}\vec{Q}^\T) \vec{A} (\vec{I} - \vec{Q}\vec{Q}^\T) 
+ (\vec{I} - \vec{Q}\vec{Q}^\T) \vec{A} (\vec{I} + \vec{Q}\vec{Q}^\T)  \Big].
\end{equation}
Introduce matrices
\begin{equation}\label{eqn:ZfHZ_def}
\vec{Z} := (\vec{I} - \vec{Q}\vec{Q}^\T) \vec{Y}
,\qquad
\vec{W} := (\vec{I} + \vec{Q}\vec{Q}^\T) \vec{Y}.
\end{equation}
Then, in place of \cref{eqn:diff_hard}, we can use the mathematically equivalent expression
\begin{equation}
    \frac{1}{2} \Big[ \vec{W}^\T \vec{A} \vec{Z} + \vec{Z}^\T \vec{A} \vec{W} \Big].
\end{equation}
We expect this expression to be less prone to rounding errors so long as $\vec{W}^\T \vec{A}\vec{Z}$ is far from skew-symmetric. 

A common way to efficiently obtain an orthonormal matrix $\vec{Q}$ for which the error $\|\vec{A} - \vec{Q}\vec{Q}^\T\vec{A}\vec{Q}\vec{Q}^\T\|_\F^2$ is small is to take $\vec{Q} = \operatorname{orth}(\vec{A}\vec{\Omega})$, where $\vec{\Omega}$ is a $\dt \times k$ random matrix with standard normal entries.
The resulting approximation is commonly called the randomized SVD and requires $k$ matrix-vector products with $\vec{A}$.
When $\vec{A}$ has a quickly decaying eigenvalues, the resulting $\vec{Q}$ results in an approximation nearly as good as the exact top eigenspace.
For matrices with slower decaying eigenvalues, there are more complicated algorithms \cite{halko_martinsson_tropp_11,tropp_webber_23}.

When $\vec{A} = \exp(-\beta \vec{H}_{\syst}) / \tr(\exp(-\beta \vec{H}_{\syst}))$, there are a number of challenges to using an approximate top subspace rather than the true top eigenspace.
First, computing matrix-vector products with $\vec{A}$ requires the use of some sort of iterative method using products with $\vec{H}_{\syst}$. 
While there are tools for this \cite{druskin_knizhnerman_89,gallopoulos_saad_92,moler_vanloan_03}, these tools are not as mature as eigensolvers.
Second, we would like to use a single matrix $\vec{Q}$ for all values of $\beta$ in some range of interest. 
The best $\vec{Q}$ will obtained by applying the randomized SVD to the matrix corresponding to the largest value of $\beta$; in fact, as $\beta\to\infty$, this will result in exactly obtaining the top subspace. 
Finally, and perhaps most subtly, when we exactly deflate the top eigenspace of $\vec{A}$, as noted in \cref{eqn:ZfHZ_bound}, the convergence of the iterative method used to compute products with $\vec{A}$ is accelerated. 
This acceleration is often significant, especially when $\beta$ is large. 
This means the cost savings of using an approximate top subspace must be compared with the additional overhead of subsequent products with $\vec{A}$.

\begin{figure}[t]
    \centering
    \includegraphics[scale=.55]{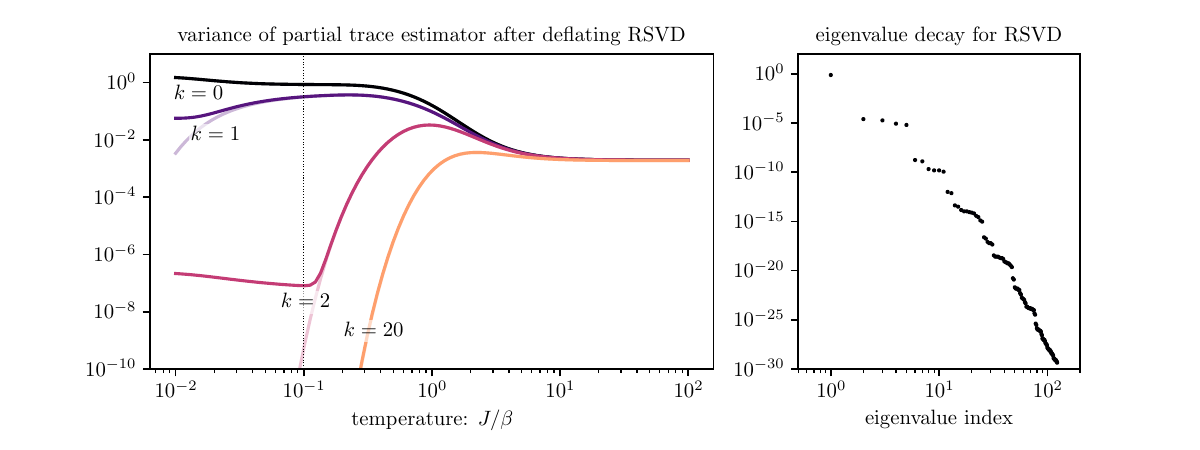}
    \caption{
    Using the same example as \cref{fig:variance_example}, we plot the approximate variance of the partial trace estimator for $\vec{A} = \exp(-\beta \vec{H}_{\syst}) / \tr(\exp(-\beta \vec{H}_{\syst}))$ when $\vec{Q} = \operatorname{orth}(\exp(-\beta_0 \vec{H}_{\syst})\vec{\Omega})$ is an approximate top subspace. We use $J/\beta_0 = 0.1$ (dotted vertical line) and sample $\vec{\Omega}\in\mathbb{R}^{\dt \times k}$ with independent standard normal entries.
    The light curves are those of \cref{fig:variance_example} and correspond to the optimal rank-$k$ subspace.
    While the approximate top subspace does reduce the variance, the variance does not go to zero in the zero temperature limit $\beta\to\infty$.
    The right plot shows the normalized singular values of $\exp(-\beta_0 \vec{H}_{\syst})$, which decay rapidly.
    }
    \label{fig:variance_example_sketch}
\end{figure}

In \cref{fig:variance_example_sketch} we plot the quantity
\begin{equation}\label{sketch_variance}
    2 \|\vec{A} - \vec{Q} \vec{Q}^\T \vec{A} \vec{Q} \vec{Q}^\T\|_\F^2
    ,\qquad 
    \vec{A} = \exp(-\beta \vec{H}_{\syst}) / \tr(\exp(-\beta \vec{H}_{\syst})).
\end{equation}
where $\vec{Q}$ is obtained by applying the randomized SVD to $\exp(-\beta_0 \vec{H}_{\syst})$, for some fixed value $\beta_0$.
This is compared to \cref{fig:variance_example}, where we plot the analogous quantity when the projection space is taken as the top subspace.
For $\beta < \beta_0$, the variance reduction essentially the same as if the exact top eigenspace were used.
For $\beta > \beta_0$ there is still some variance reduction, it does not result in a zero-variance approximation in the zero-temperature limit $\beta\to\infty$.

A deeper exploration of the tradeoffs between the cost to compute $\vec{Q}$, the quality of the variance reduction, and the costs of computing $f(\vec{H})\vec{Q}$ is beyond the scope of the present paper, but is an important topic for future work.

\section{Numerical experiments}
\label{sec:numerical}

Our experiments focus on Heisenberg spin systems in an isotropic magnetic field oriented in the positive z-direction:
\begin{equation}
    \label{eqn:hamiltonian}
    \vec{H} := \sum_{i,j=1}^{N} \left[  J^{\textup{x}}_{i,j} \bm{\sigma}^{\textup{x}}_i \bm{\sigma}^{\textup{x}}_j 
    +J^{\textup{y}}_{i,j} \bm{\sigma}^{\textup{y}}_i \bm{\sigma}^{\textup{y}}_j
    +J^{\textup{z}}_{i,j} \bm{\sigma}^{\textup{z}}_i \bm{\sigma}^{\textup{z}}_j \right] 
    + \frac{h}{2} \sum_{i=1}^{N} \bm{\sigma}_i^{\textup{z}}.
\end{equation}
Here \( \bm{\sigma}^{\textup{x}/\textup{y}/\textup{z}}_i \) is defined by
\begin{align}
    \bm{\sigma}^{\textup{x/y/z}}_i
    = \underbrace{\vec{I} \otimes \cdots \otimes \vec{I}}_{i-1\text{ terms}} 
    \otimes ~ \bm{\sigma}^{\textup{x/y/z}} \otimes 
    \underbrace{\vec{I} \otimes \cdots \otimes \vec{I}}_{N-i\text{ terms}},
\end{align}
where \( \bm{\sigma}^{\textup{x}/\textup{y}/\textup{z}} \) are the Pauli spin-$\frac{1}{2}$ matrices
\begin{align}
    \bm{\sigma}^{\textup{x}} = \begin{bmatrix} 0 & 1 \\ 1 & 0 \end{bmatrix}
    && 
    \bm{\sigma}^{\textup{y}} = \begin{bmatrix} 0 & -\ii \\ \ii & 0 \end{bmatrix}
    && 
    \bm{\sigma}^{\textup{z}} = \begin{bmatrix} 1 & 0 \\ 0 & -1 \end{bmatrix}.
\end{align}
We remark that while $\bm{\sigma}^{\textup{y}}$ is Hermitian, $\bm{\sigma}^{\textup{y}}_i\bm{\sigma}^{\textup{y}}_j$ (and thus $\vec{H}$) is \emph{real symmetric}.

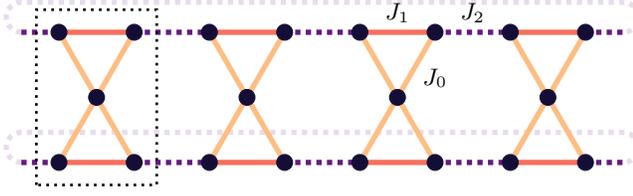
\begin{figure}[h]
\input{kagome}
    \caption{Kagome-strip chain with $N=20$ sites and periodic boundary conditions. Subsystem (s) is encircled.}
    \label{fig:kagome_graph}
\end{figure}

\subsection{Experimental setup}
Our experiments are implemented in Python using double precision arithmetic. 
We use scipy's sparse library to represent Hamiltonians and scipy's \texttt{sparse.linalg.eigsh} to compute the top eigenvectors. 
The latter is a wrapper for ARPACK's \texttt{dsaupd}, which is an implementation of the Implicitly Restarted Lanczos Method, and computes the eigenvectors to machine precision \cite{lehoucq_sorensen_yang_98}.

We set the number of Lanczos iterations $t$ so that products involving $\exp(-\beta \vec{H})$ are computed to a relative error roughly $10^{-10}$. 
The statistical noise from the random samples is much larger, so computing the matrix functions any additional accuracy does not impact the results in any noticeable way.
Re-orthogonalization is not used, but we do explicitly orthogonalize against the deflated subspace.

Code used to generate the data plotted in the figures is available at \texttt{\url{https://github.com/tchen-research/faster_partial_trace}}.

\subsection{Quantities of interest}

\subsubsection{Entanglement spectrum and von Neumann entropy}

The von Neumann entropy is an information theoretic measure of the entropy of a quantum system and can be viewed as a measure of how far a quantum state is from pure.
Thus, the von Neumann entropy of subsystem (s) provides information about the entanglement between subsystems (s) and (b).

The von Neumann entropy of subsystem (s) is defined by the formula
\begin{equation}
S = S(\beta,h) := - \tr\big(\bm{\rho}^* \ln(\bm{\rho}^*)\big),
\end{equation}
where $\bm{\rho}^*$ is the reduced density matrix defined in \cref{eqn:rho_HMF}.

The set of eigenvalues of $-\ln(\bm{\rho}^*)$ is sometimes referred to as the entanglement spectrum, and provides a more complete picture of the entanglement of subsystems (s) and (b) than the von Neumann entropy \cite{li_haldane_08}. Up to a scaling factor $1/\beta$, $-\ln(\bm{\rho}^*)$ is the same as the Hamiltonian of the mean force, and important quantity in equilibrium thermodynamics \cite{talkner_hanggi_20}.

\subsubsection{Ergotropy of quantum batteries}
Quantum batteries use quantum systems to store energy, and offer the potential for faster charging and higher efficiency than classical batteries \cite{campaioli_pollock_vinjanampathy_18}.
We consider the setup of \cite{barra_hovhannisyan_imarato_22} in which the battery consists of the spins in subsystem (s) and charges the spins in subsystem (b).
Once the total system (battery+charger) is in thermal equilibrium, the battery is instantaneously disconnected from the charger and is therefore in state $\bm{\rho}^*$ with internal energy $\tr(\vec{H}_{\syss}\bm{\rho}^*)$.
We can extract energy from the battery by evolution with a unitary $\vec{U}$, which brings us to state $\vec{U}\bm{\rho}^* \vec{U}^\T$ with internal energy $\tr(\vec{H}_{\syss}  \vec{U}\bm{\rho}^* \vec{U}^\T)$.
The ergotropy \cite{allahverdyan_balian_nieuwenhuizen_04,campaioli_gherardini_quach_polini_andolina_23} $\mathcal{E} = \mathcal{E}(\beta,h)$, is defined as the total possible energy which could be extracted from the battery:
\begin{equation}
\mathcal{E} = \mathcal{E}(\beta,h) := \max_{\vec{U}^\cT\vec{U} = \vec{I}} \bigg( \tr\big(\vec{H}_{\syss} \bm{\rho}^*\big) -  \tr\big(\vec{H}_{\syss}  \vec{U}\bm{\rho}^* \vec{U}^\T\big) \bigg).
\end{equation}
The unitary $\vec{U}$ minimizing $\tr\big(\vec{H}_{\syss}  \vec{U}\bm{\rho}^* \vec{U}^\T\big)$ can be obtained explicitly \cite{allahverdyan_balian_nieuwenhuizen_04}.   
Specifically, if $\vec{H}_{\syss}$ and $\vec{\bm{\rho}}^*$ are diagonalized (with eigenvalues in non-increasing order) as $\vec{H}_{\syss} = \vec{Q}_{\syss} \vec{\Lambda}_{\syss} \vec{Q}_{\syss}^\T$ and $\bm{\rho}^* = \vec{Q}_{\rho} \vec{\Lambda}_{\rho} \vec{Q}_{\rho}^\T$, then 
\begin{equation}
\vec{U} =  \vec{Q}_{\syss} \vec{P} \vec{Q}_{\rho}^\T,
\end{equation}
where $\vec{P}$ is the reversal permutation matrix (identity with columns reversed).

\subsection{Kagome-strip chain}

In this experiment, we consider Kagome-strip chain systems \cite[etc.]{azaria_hooley_lecheminant_lhuillier_tsvelik_98,white_singh_00,morita_sugimoto_sota_tohyama_18} as show in \cref{fig:kagome_graph}.
We take subsystem (s) to be the 5 spins indicated in \cref{fig:kagome_graph}.

We choose several values of $J_2$ and fix $J_1 = J_3=J$.
For each choice of $J_2$, we consider a range of $h$ and $\beta$ for each system.
To determine which values of $h$ to run our algorithm at, we use bisection on the von Neumann entropy of the ground state to determine intervals where the von Neumann entropy appears constant.
We then run our algorithm at values of $h$ corresponding to Chebyshev nodes shifted and scaled to each interval.
Throughout, we use $k=25$ and $m=5$.

\begin{figure}[htb]
    \centering
    \includegraphics[scale=0.6]{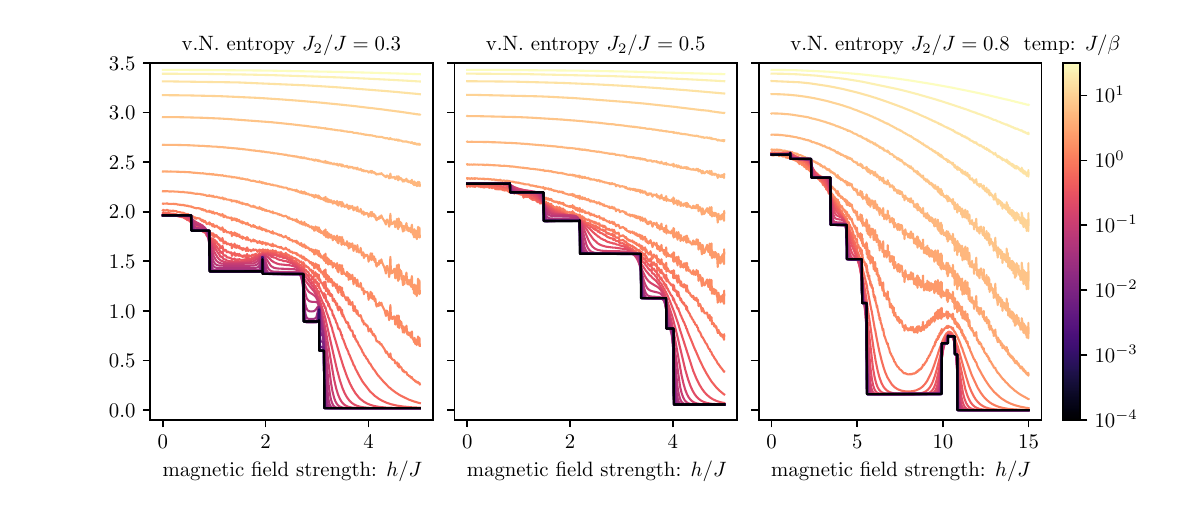}
    \caption{Von Neumann entropy for Kagome-strip chain at varying values of $J_2$, $h$, and $\beta$.}
    \label{fig:kagome_entropy}
\end{figure}

\begin{figure}[htb]
    \centering
    \includegraphics[scale=0.6]{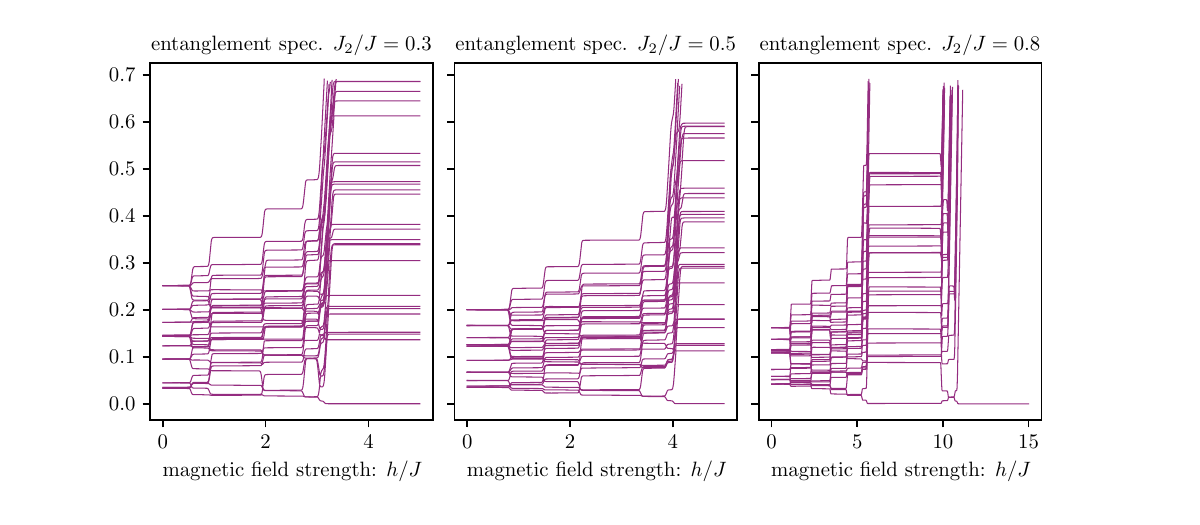}
    \caption{Entanglement spectrum for Kagome-strip chain at varying values of $J_2$ and $h$ at fixed temperature $J/\beta = 2\times 10^{-2}$.} \label{fig:entanglement_spec_KSC}
\end{figure}

We also consider the entanglement spectrum at a fixed value of $\beta$. 
This is illustrated in \cref{fig:entanglement_spec_KSC}.

\subsection{Long range spin chain}
\label{sec:spin_change_lr}

We now consider the XX spin chain with long-range power-law interactions
\begin{equation}     
J^{\textup{x}}_{i,j} =  J^{\textup{y}}_{i,j} = |i-j|^{-\alpha}
     ,\qquad  J^{\textup{z}}_{i,j} = 0.
\end{equation}
Here we take subsystem (s) to be the first 2 spins, and subsystem (b) to be the remaining spins.
We remark that in the case $\alpha = \infty$, this system is exactly solvable, and we use this to verify the accuracy of our algorithm in \cref{sec:solvable_chain}.

Within this framework, we set $N=16$ and vary the parameters $\alpha$, $h$, and $\beta$.
We use the same bisection-based approach to determine suitable values of $h$ to run our algorithm at. 
For each value of $\alpha$ and $h$, we run our algorithm with $k=25$ and $m=5$ and again use enough Lanczos iterations to accurately compute the matrix functions for each value of $\beta$.

\begin{figure}[tb]
    \centering
    \includegraphics[scale=.55]{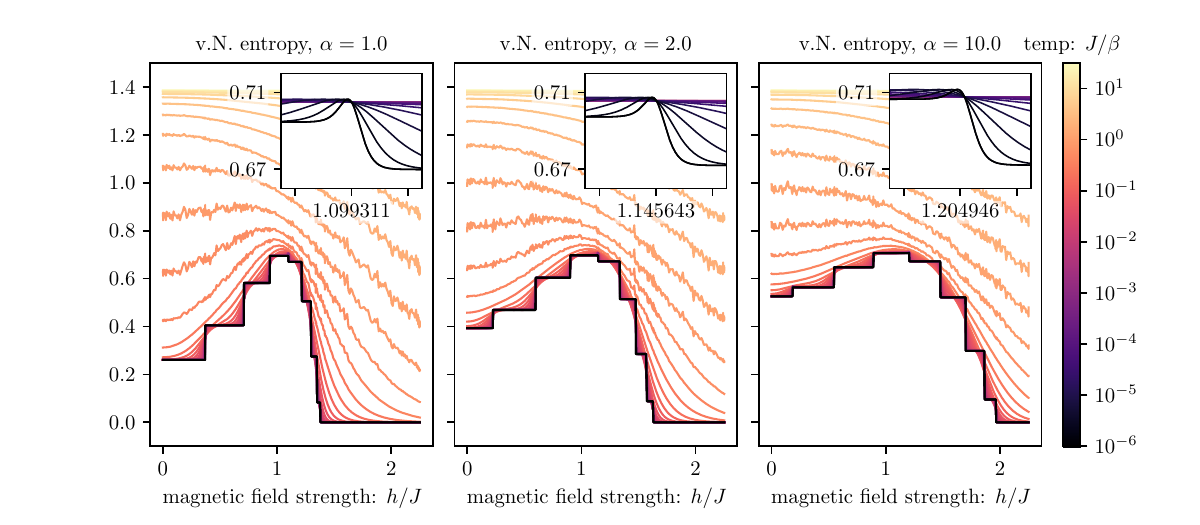}
    \caption{Von Neumann entropy for long-range spin chain with $N=16$ and the system taken as the first two spins at varying values of $\alpha$, $h$, and $\beta$. Inset figure shows the transition to the right of the top ``plateau''.}
    \label{fig:vN_phase_alpha}
\end{figure}

\begin{figure}[b!]
    \centering
    \includegraphics[scale=.55]{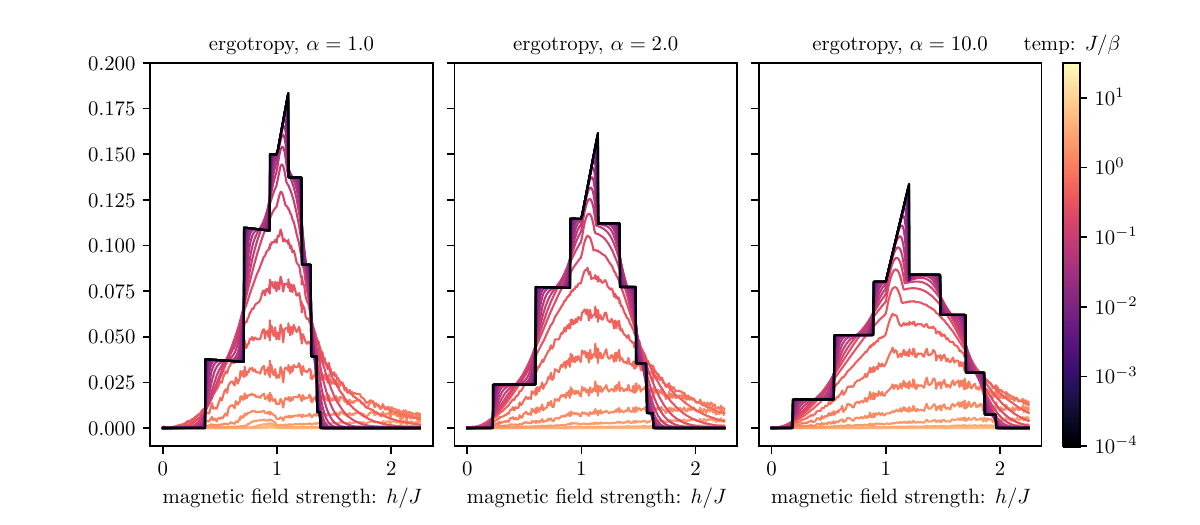}
    \caption{Ergotropy for long-range spin chain with $N=16$ and the system taken as the first two spins at varying values of $\alpha$, $h$, and $\beta$. Observe the non-piecewise continuous behavior to the right of 1, despite the von Neumann entropy appearing constant in this region.}
    \label{fig:ergotropy_lr_combined}
\end{figure}
In \cref{fig:vN_phase_alpha} we visualize the von Neumann entropy of subsystem (s) for several values of $\alpha$.
We observe that at zero temperature, the von Neumann entropy appears piecewise constant.
In the solvable model, the steps in the zero temperature von Neumann entropy correspond to values of $h/J$ for which a fermionic eigenmode vanishes \cite[(69)]{campisi_zueco_talkner_10}.
The inset panel of \cref{fig:vN_phase_alpha} shows a plot zoomed in to the right edge of the top ``plateau'' and illustrates that the von Neumann entropy of the ground state changes continuously in this region.


In \cref{fig:ergotropy_lr_combined} we show the ergotropy of subsysem (s) for several values of $\alpha$.
We use the same values of $h$ used for the von Neumann entropy.
While the sharp increases appear to happen in the same places as the von Neumann entropy, the regions between jumps appear linear rather than constant.
In addition, there is an apparent discontinuity in the derivative of the $\beta = \infty$ curve at $h/J = 1$, which does not appear in the von Neumann entropy.

\begin{figure}[htb]
    \centering
    \includegraphics[scale=.55]{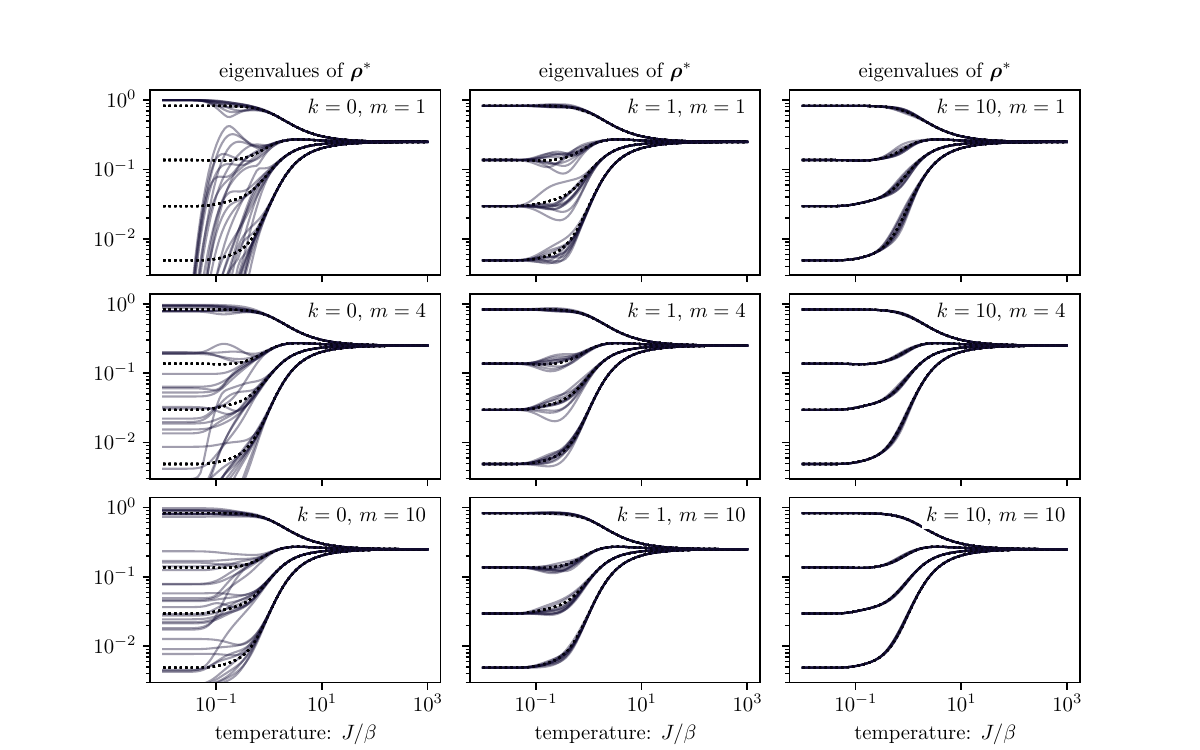}
    \caption{Comparison of $k=0$ (equivalent to \cite{chen_cheng_22}), $k=1$ and $k=10$ for several values of $m$ when our algorithm is run on the solvable model with $N=16$ and the system taken as the first two spins 
    Lines correspond to repeated runs of our algorithm. 
    \emph{Takeaway}: As $k$ and $m$ increase, the variance decreases. However, while the variance decreases linearly with $m$, it may decrease more quickly with $k$.}
    \label{fig:solvable_compare}
\end{figure}

\subsection{Validation on the solvable XX spin chain}
\label{sec:solvable_chain}

In the special case that
\begin{equation}
    \label{eqn:solvable}
     J^{\textup{x}}_{i,j} =  J^{\textup{y}}_{i,j} 
     = \begin{cases} J & |i-j| = 1 \\ 0 & |i-j|\neq 1 \end{cases}
     ,\qquad  J^{\textup{z}}_{i,j} = 0,
\end{equation}
the system \cref{eqn:hamiltonian} is exactly solvable via the ``Bethe ansatz'' \cite{karabach_muller_gould_tobochnik_97}.
That is, it can be diagonalized analytically.
In addition, expressions for the partial trace of the first 2 spins have been obtained \cite{campisi_zueco_talkner_10}; see also \cite[Appendix C]{chen_cheng_22}.
This allows us to test our algorithm against a known solution for problem sizes where exact diagonalization is intractable.

\subsubsection{Variance study}

We begin by studying the variability of the output of our algorithm.
We use $N=18$ and $h/J = 0.3$ and run the algorithm at varying values of $k$ and $m$ and compute the eigenvalues of $\vec{\rho}^*$ for a range of $\beta$.
For each choice of $k$ and $m$, we repeatedly and independently run our algorithm 10 times.
This gives some indication of the variance in the algorithm. 
In all cases, we use enough Lanczos iterations to accurately compute the matrix functions for each value of $\beta$.
\Cref{fig:solvable_compare} shows the results of this experiment.
Note that the $k=0$ plots in \cref{fig:solvable_compare} correspond to the algorithm from \cite{chen_cheng_22} which does not use deflation; see \cite[Figure 1]{chen_cheng_22}.

\subsubsection{Jackknife variance estimates}

\begin{figure}[t]
    \centering
    \includegraphics[scale=.55]{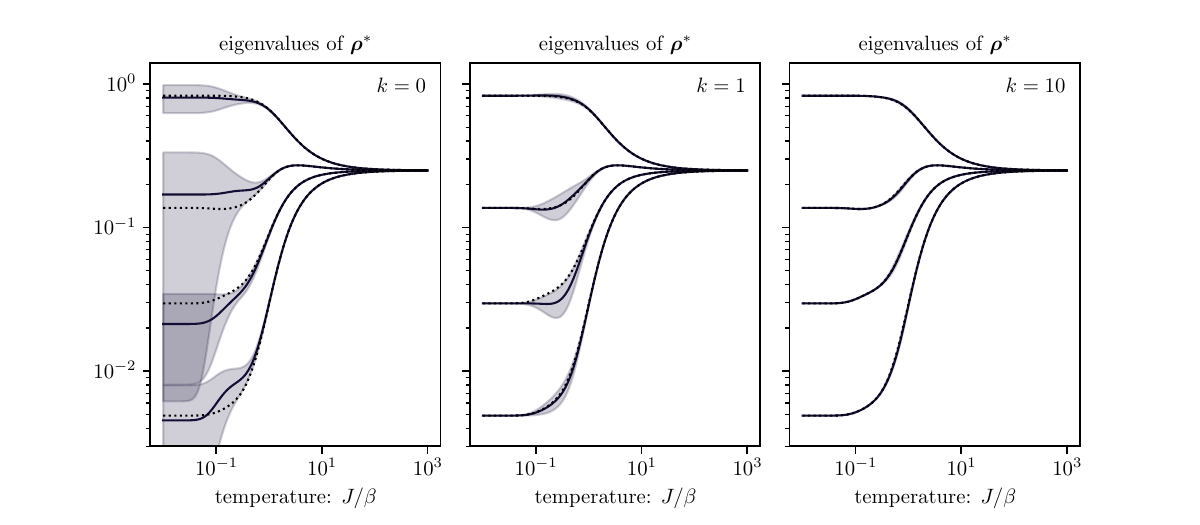}
    \caption{Leave-one-out estimator for standard error obtained from a single run of the algorithm (same setup as \cref{fig:solvable_compare}) with $m=10$.
    \emph{Takeaway}: Using just the information from a single run of the algorithm, we can get reasonable estimates for the variability of the output.}
    \label{fig:solvable_compare_est_errs}
\end{figure}

Generating plots like \cref{fig:solvable_compare} is not practical, as it requires multiple runs of the algorithm for each parameter setting.
However, since our estimator \cref{eqn:residual_low-rank} involves averaging $m$ iid samples of an unbiased random matrix, we can use a Jackknife (leave-one-out) estimator for the variance. 
We visualize the error estimated by the Jackknife method for a single run of the algorithm \cref{fig:solvable_compare_est_errs}.
Here, the estimated standard error seems to align well with the true error of the algorithm.

\subsubsection{Von Neumann Entropy}

\begin{figure}[htb]
    \centering
    \includegraphics[scale=.55]{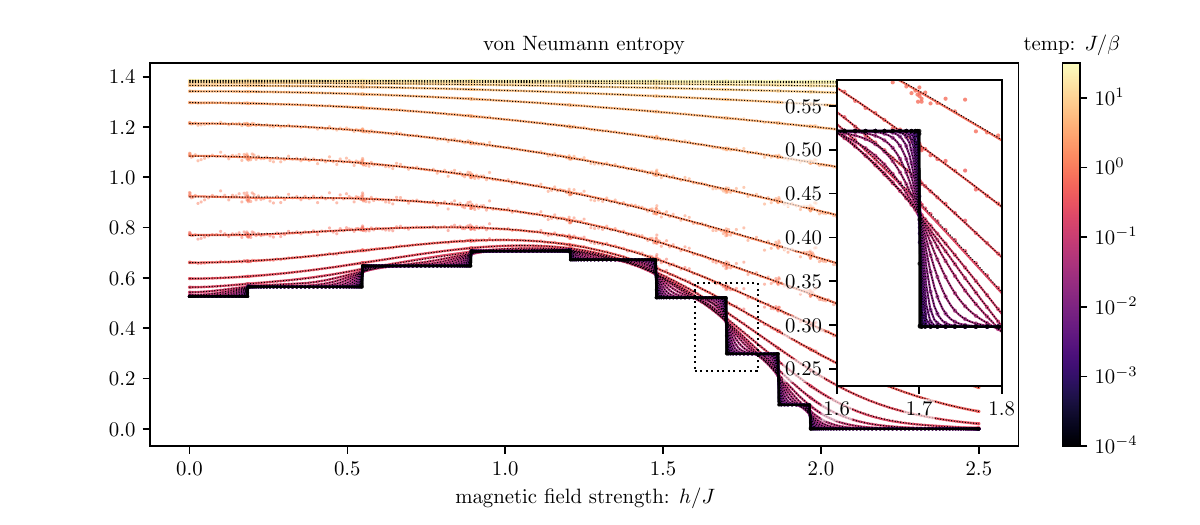}
    \caption{Approximation of von Neumann entropy on the solvable model with $N=16$ and the system taken as the first two spins. 
    Smoothing is used to reduce noise at moderate temperatures (raw data shown as dots).
    We use $k=25$ eigenvectors for deflation and $m=5$ copies of the stochastic trace estimator. 
    Compared to \cite[Figure 3]{chen_cheng_22} $(k=0,m=400)$, our approach is both computationally cheaper and results in a much more accurate approximation.
    \emph{Takeaway}: Deflation allows the low-temperature parts of the curve to be obtained very accurately. While the algorithm's output has some noise at intermediate temperature, this is easily overcome by smoothing, which naturally averages the randomness of nearby points together.}
    \label{fig:vN_entropy}
\end{figure}

In this experiment, we set $N = 16$ and $N_{\syss} = 2$ and vary the magnetic field strength $h/J$.
\Cref{fig:vN_entropy} shows the exact von Neumann entropy as a function of $h/J$ for a range of $\beta$ as well as the values computed by our algorithm with the parameters $k=25$ and $m=5$.
In the right panel of \cref{fig:vN_entropy}, we show a cropped version of the left panel. 
In this plot, the low-variance behavior of the algorithm at low but nonzero temperatures is clearly visible. 
This is in sharp contrast to \cite{chen_cheng_22}, in which high variance is observed at low temperatures, even with $m=400$.

For high temperatures $\beta/J <  5$, we simply fit a degree 10 polynomial with least squares. 
While the raw data deviates considerably from the true von Neumann entropy, the least squares fit seems extremely accurate. 
This suggests that the bias of the algorithm's output is fairly small.
For low temperature $\beta/J \in (5,500)$, we use cubic splines. 
Here the algorithm's output has little noise, and the smoothing is mainly to interpolate the data to values of $h/J$ which we did not run the algorithm on.
Finally, for very low temperature $\beta/J > 500$, we do not do any smoothing.

\section{Conclusion and Outlook}

We have incorporated projection as a means of variance reduction in typicality estimators for the partial traces. 
Since the partial trace does not satisfy the cyclic property, there is a potential for cancellation if the partial trace estimator cannot be computed exactly.
To avoid this, we use deflation (projection of top eigenspace) in combination with the block-Lanczos algorithm and explicitly orthogonalize against the projected subspace.
Our approach significantly reduces the runtime required to obtain accurate estimates, often by several orders of magnitude.

In the future, we would like to take further advantage of the structure of the systems we are dealing with.
For example, in many situations, one would like to compute a quantity of interest over a range of parameter values (e.g. magnetic field strength, coupling strength, etc.). 
Presently, we apply our algorithm independently at each parameter setting. 
This leaves open the potential for a faster algorithm which can take advantage of the fact that many of the quantities may not change significantly.

\bibliographystyle{siamplain}
\bibliography{refs}
\end{document}

%% file: kagome.tex
    \centering
\def\dx{1}
\def\dy{sqrt(3/4)}
\def\offx{cos(60)}
\def\loopy{.2}
\begin{tikzpicture}

    \draw[line width=2pt,dotted,c1,opacity=.15] ({0*\dx+\offx},{2*\dy}) arc (270:90:\loopy)  -- ({8*\dx+1*\offx},{2*\dy+\loopy*2}) arc (90:-90:\loopy);
    \draw[line width=2pt,dotted,c1,opacity=.15] ({0*\dx+\offx},{0}) arc (270:90:\loopy)  -- ({8*\dx+1*\offx},{0+\loopy*2}) arc (90:-90:\loopy);
    
    \foreach \i in {2,4,...,8} {
        \draw[line width=2pt,c3] ({\i*\dx-\offx},{\dy}) -- ({\i*\dx},0);
        \draw[line width=2pt,c3] ({\i*\dx-\offx},{\dy}) -- ({\i*\dx-2*\offx},0);
        \draw[line width=2pt,c3] ({\i*\dx-\offx},{\dy}) -- ({\i*\dx},{2*\dy});
        \draw[line width=2pt,c3] ({\i*\dx-\offx},{\dy}) -- ({\i*\dx-2*\offx},{2*\dy});

        \draw[line width=2pt,c2] ({\i*\dx-2*\offx},{2*\dy}) -- ({\i*\dx},{2*\dy});
        \draw[line width=2pt,c2] ({\i*\dx-2*\offx},{0}) -- ({\i*\dx},{0});
    }
    \foreach \i in {2,4,6} {
        \draw[line width=2pt,dotted,c1] ({\i*\dx},{2*\dy}) -- ({\i*\dx+2*\offx},{2*\dy});
        \draw[line width=2pt,dotted,c1] ({\i*\dx},{0}) -- ({\i*\dx+2*\offx},{0});
    }
    \draw[line width=2pt,dotted,c1] ({8*\dx},{2*\dy}) -- ({8*\dx+1*\offx},{2*\dy});
    \draw[line width=2pt,dotted,c1] ({8*\dx},{0}) -- ({8*\dx+1*\offx},{0});
    \draw[line width=2pt,dotted,c1] ({0*\dx+\offx},{2*\dy}) -- ({0*\dx+2*\offx},{2*\dy});
    \draw[line width=2pt,dotted,c1] ({0*\dx+\offx},{0}) -- ({0*\dx+2*\offx},{0});

    \foreach \i in {2,4,...,8} {
        \filldraw[c0] ({\i*\dx-\offx},{\dy}) circle (3pt);
    }
    \foreach \i in {1,2,...,8} {
        \filldraw[c0] ({\i*\dx},{0}) circle (3pt); 
        \filldraw[c0] ({\i*\dx},{2*\dy}) circle (3pt); 
    }
    \draw[line width=1,dotted] ({2*\dx+\dx-.3},-.3) rectangle ({2*\dx+2*\dx+.3},{2*\dy+.3});
    \node[] at ({6*\dx},{1.3*\dy}) {\small$J_0$};
    \node[] at ({5.5*\dx},{2.3*\dy}) {\small$J_1$};
    \node[] at ({6.5*\dx},{2.3*\dy}) {\small$J_2$};
\end{tikzpicture}